\documentclass[a4paper,12pt]{amsart}
\usepackage{amsmath}
\usepackage{amsthm}
\usepackage{amscd}
\usepackage{amssymb}
\usepackage{amsfonts}
\usepackage{latexsym}
\usepackage[mathscr]{eucal}

\newtheorem{thm}{Theorem}[section]
\newtheorem{prop}[thm]{Proposition}
\newtheorem{lemma}[thm]{Lemma}
\newtheorem{cor}[thm]{Corollary}

\def\supp{\mathop{\rm {supp}}\nolimits}
\def\graph{\mathop{\rm {graph}}\nolimits}
\def\length{\mathop{\rm {length}}\nolimits}

\begin{document}
\title{$C^*$-algebras associated with self-similar sets}
\author{Tsuyoshi Kajiwara}
\address[Tsuyoshi Kajiwara]{Department of Environmental and 
Mathematical Sciences, 
Okayama University, Tsushima, 700-8530,  Japan}      

\author{Yasuo Watatani}
\address[Yasuo Watatani]{Department of Mathematical Sciences, 
Kyushu University, Hakozaki, 
Fukuoka, 812-8581,  Japan}
\maketitle
\begin{abstract}
Let $\gamma = (\gamma_1,\dots,\gamma_N)$, $N \geq 2$,  
be a system of  
proper contractions on a complete metric space.  
Then there exists a unique 
self-similar non-empty compact subset $K$.  We consider 
the union ${\mathcal G} =  
\cup _{i=1}^N \{(x,y) \in K^2 ; x = \gamma _i(y)\}$
of the cographs of $\gamma _i$. Then 
$X = C({\mathcal G})$ is a Hilbert bimodule over $A = C(K)$.
We associate a $C^*$-algebra 
${\mathcal O}_{\gamma}(K)$ with them as a 
Cuntz-Pimsner algebra ${\mathcal O}_X$. 
We show that if a system of proper contractions satisfies the 
open set condition in $K$, then the $C^*$-algebra 
${\mathcal O}_{\gamma}(K)$   
is simple and purely infinite, which is not isomorphic to 
a Cuntz algebra in general. 

\end{abstract}

\section{Introduction} The study of the self-similar set 
constructed from iterations of proper contractions has deep 
interaction with many areas of mathematics. The 
theory of $C^*$-algebras seems to be one of them. For example 
Bratelli-Jorgensen \cite{BJ} considered a relation among
representations of the Cuntz algebra \cite{C}, wavelet theory and  
the iterated function systems. See also \cite{MSW}. 
In this paper 
we shall give a new construction of a $C^*$-algebra 
associated with a system of proper 
contractions on a self-similar set.  
The algebra is not a Cuntz algebra in general and its K-theory 
is closely related with the failure of the injectivity 
of the coding by the full shift. In the case that the contractions 
are branches of the inverse of a certain map $h$, 
its K-theory is related with 
the structure of the branched points (critical points) of the map $h$.    

Let $\gamma = (\gamma_1,\dots,\gamma_N)$, $N \geq 2$,  be a system of 
proper contractions 
on a complete metric space $\Omega$.  Then there exists a unique 
compact non-empty subset $K \subset \Omega$ satisfying 
the self-similar condition such that 
$K  = \cup _i \gamma_i(K)$. In the paper we usually 
foget an ambient space $\Omega$ and regard 
each $\gamma_i$ is a map on $K$. The subset 
$\{(x,y) \in K^2 ; x = \gamma _i(y)\}$ of $K^2$ is called 
the {\it cograph} of $\gamma _i$.  Define  
${\mathcal G} =  
\cup _{i=1}^N \{(x,y) \in K^2 ; x = \gamma _i(y)\}$
be the union of the cographs of $\gamma _i$. If the contractions 
are the continuous branches of the inverse of a certain map 
$h: K \rightarrow K$, then ${\mathcal G}$ is exactly the graph 
of $h$. Let $A = C(K)$ be the algebra of continuous functions 
on the self-similar set $K$.
Define an endomorphism 
$\beta _i : A \rightarrow A$ by 
$(\beta _i (a))(y) = a(\gamma _i (y))$
for $a \in A$, $y \in K$. Let $C^*(A,\beta_1,\dots, \beta_N)$ 
be the universal $C^*$-algebra generated by $A$ and the Cuntz 
algebra ${\mathcal O}_N = C^*(S_1,\dots, S_N)$ with the 
commutation relations $aS_i = S_i\beta_i(a)$ for $a \in A$ and 
$i = 1,\dots, N$. Since each $\gamma_i$ is a proper contraction, 
$C^*(A,\beta_1,\dots, \beta_N)$ turns out to be isomorphic to 
the Cuntz algebra ${\mathcal O}_N$ itself as we considered  
in \cite{PWY}. 
The problem of this construction is that we forgot to pay 
attention to the "branched points" and used  
the simpler object, i.e.,  
the {\it disjoint} union of the the cographs of $\gamma _i$ 
instead of the usual union ${\mathcal G} \subset K^2$ of the 
cographs of $\gamma _i$. In the language of bimodule, the disjoint 
union corresponds to a direct sum $\oplus _i  \ _{\beta_i} A$ 
of bimodules and the union ${\mathcal G}$ 
correponds to a bimodule $X= C({\mathcal G})$ which is 
embedded as a submodule of $\oplus _i  \ _{\gamma_i} A$.  
Under the consideration, in the paper we shall give a new construction  
of a $C^*$-algebra ${\mathcal O}_{\gamma}(K)$ associated 
with a system $\gamma = (\gamma_1,\dots,\gamma_N)$ on K as a 
Cuntz-Pimsner algebra \cite{Pi} of a Hilbert bimodule 
$X = C({\mathcal G})$ over the algebra $A = C(K)$. Then the 
$C^*$-algebra ${\mathcal O}_{\gamma}(K)$ is not 
a Cuntz algebra in general. In fact the $K_0$-group can have a 
torsion free element. If two systems of contractions are topologically 
conjugate, then associated algebras are isomorphic.  
In a recent paper \cite{KW} we introduced the $C^*$-algebra 
${\mathcal O}_R(J_R)$ associated with a rational function $R$ on  
its Julia set $J_R$ using bimodules, after pioneering 
works Deanonu \cite{D} and 
Deaconu and Muhly \cite{DM} on $C^*$-algebras associated with 
branched coverings using 
a groupoid approach of Renault \cite{R}.  
If the inverse of the rational function $R$ 
has continuous branches $\gamma_1,\dots,\gamma_N$ on $J_R$, 
then our construction is arranged so that $C^*$-algebra 
${\mathcal O}_R(J_R)$ is isomorphic to the $C^*$-algebra 
${\mathcal O}_{\gamma}(J_R)$, because 
$\graph R = \cup_i cograph \ \gamma_i$.  We note that there 
exists an example of a rational functon $R$ that the Julia set $J_R$ 
is homeomorphic to the Sierpinski gasket \cite{Kam}, \cite{U}.  
We show that if a system of proper contractions satisfies the 
open set condition in $K$, then the $C^*$-algebra 
${\mathcal O}_{\gamma}(K)$   
is simple and purely infinite.  

We remark that there exists an analogy such that 

\ Klein group : limit set : crossed product \\
=  rational function $R$ : Julia set $J_R$ : ${\mathcal O}_R(J_R)$ \\
=  system $\gamma$ of contractions : self-similar set $K$ :
 ${\mathcal O}_{\gamma}(K)$.

See also Anatharaman-Delaroch \cite{A} and Laca-Spielberg \cite{LS}, 
and Kumujian \cite{Ku} for constructions of purely infinite, simple 
$C^*$-algebras.

We also note that the $C^*$-algebra 
${\mathcal O}_{\gamma}(K)$ is related with 
graph $C^*$-algebras \cite{KPRR} and their generalizaion for  
topological relations by Brenken \cite{Br}, topological graphs by 
Katsura \cite{Ka1}, \cite{Ka2} and topological quivers by 
Muhly and Solel \cite{MS} and by Muhly and Tomforde \cite{MT}.   

After almost completing our work,  we have found a preprint 
\cite{MT} by Muhly and Tomforde. The simplicity and some other properties of 
the $C^*$-algebra ${\mathcal O}_{\gamma}(K)$ except purely infiniteness 
follows from their general work. We give a direct proof of simplicity and 
purely infiniteness. 

We have also learned that Nekrashevych 
have introduced interesting $C^*$-algebras associated 
with graph-directed iterated 
function systems in a survey paper \cite{BGN}. In particular case, 
if maps are proper contractions,  Ionesc \cite{I}
has shown that the $C^*$-algebra is in fact isomorphic to the Cunz-Kreiger 
algebras associated to the underlying finite graph as a generarization 
of \cite{PWY}. He also obtains a "no go" theorem. If the set of vertices 
is a singleton, then it is a usual iterated function system and their 
construction corresposnds to the disjoint union of the the cographs of 
$\gamma _i$. 

The authors express our thanks to M. Ionescu for pointing  out an error in 
the first verson of the present article. 

\section{Self-similar sets and Hilbert bimodules}
   Let $(\Omega,d)$ be a (separable) complete metric space $\Omega$ with a 
metric $d$.   
A map $\gamma$ on $\Omega$ is called a contraction if its Lipschitz constant 
$Lip(\gamma ) \leq 1$, that is, 
\[
Lip(\gamma ) := \sup _{x \not= y} 
\frac{d(\gamma (x), \gamma (y))}{d(x,y)}  \leq 1.
\]  
We call contractions $\{\gamma_j:j=1,2,..,N \}$   on $\Omega$ are 
{\it proper} if there exist positive constants $\{c_i\}$ and $\{c'_i\}$
with $(0 < c_i \leq c'_i < 1)$
satisfying the condition:
\[
c_i d(x,y)\leq d(\gamma_i(x),\gamma_i(y))\leq c'_i d(x,y)
\]
for any  $\ x,y \in X \quad  i=1,2,..,N$. 

We say that a non-empty compact set $K \subset \Omega$ is {\it self-similar}
(in a weak sense) 
with respect to a system $\gamma = (\gamma_1,\dots , \gamma_N)$ 
if $K$ is a finite union 
of its small copies $\gamma_j(K)$, that is,
$$
K=\cup_{i=1}^{N} \gamma_{i}(K).
$$

If the  contractons are  proper, then there exists a self-similar set 
$K \subset X$ uniquely. We say that $K$ is the self-similar set 
with respect to a system $\gamma = (\gamma_1,\dots , \gamma_N)$ and 
the self-similar set $K$ is also denoted by 
$K(\gamma) = K(\gamma_{1}, \dots , \gamma_{N})$. See \cite{F} and 
\cite{Kig} more on fractal sets.

Let $N$ be a natural number with $N \geq 2$.  For a natural number 
$m$, we define $W_m := \{1, \dots , N\}^m$ and 
$w = (w_1,\dots,w_m) \in W_m$ is called a word of length $m$ with 
symbols $\{1,2,\dots ,N\}$.  We set $W = \cup _{m \geq 1} W_m$  
and denote the length of a word $w$ by $\ell (w)$.  

The full $N$-shift space $\{1,2,...,N\}^{\mathbb N}$ is the  space 
of one-sided sequences $x = (x_n)_{n \in {\mathbb N}}$ of symbols
 $\{1,2,...,N\}$. 
We define a metric $d$ on $\{1,2,...,N\}^{\mathbb N}$ by 
\[
d(x,y) =\sum _n \frac{1}{2^n} (1 - \delta _{x_n,y_n}) .  
\] 
Then $\{1,2,...,N\}^{\mathbb N}$ is a compact metric space.  
Define a system $\{\sigma_j:j=1,2,..,N \}$ of $N$ contractions  on 
$\{1,2,...,N\}^{\mathbb N}$ by 
$$
\sigma_j(x_1,x_2,..., ) = (j,x_1,x_2,..., )
$$
Then each $\sigma_j$ is a proper contraction with the Lipschitz 
constant $Lip(\sigma _j) = \frac{1}{2}$.  
The self-similar set  $K(\sigma_{1},\sigma_{2},..,\sigma_{N}) = 
\{1,2,...,N\}^{\mathbb N}$.  

Moreover for $w = (w_1,\dots,w_m) \in W_m$, 
let  $\gamma _w = \gamma _{w_1} \circ \dots \circ \gamma_{w_m}$ and 
$K_w = \gamma _w(K)$. Then for any one-sided sequence 
$x = (x_n)_{n \in {\mathbb N}} \in \{1,2,...,N\}^{\mathbb N}$, 
$\cap _{m \geq 1} K_{(x_1, \dots ,x_m)}$ contains only one point 
$\pi (x)$.  Therefore we can define a map 
$\pi : \{1,2,...,N\}^{\mathbb N} \rightarrow K$ by 
$\{\pi (x)\} =  \cap _{m \geq 1} K_{(x_1, \dots ,x_m)}$.  
Since $\pi (\{1,2,...,N\}^{\mathbb N})$ is also a  self-similar 
set, we have $\pi (\{1,2,...,N\}^{\mathbb N}) = K$.   
Thus $\pi$ is a continuous onto map satisfying 
$\pi \circ \sigma _i = \gamma _i \circ \pi$ for $i = 1, \dots , N$. 
Moreover,  
for any $y \in K$ and any neighbourhood $U_y$ of $y$ there exists 
$n \in {\mathbb N}$ and $w \in W_n$ such that 
\[
y \in \gamma _w(K) \subset U_y .
\]

In the note  we usually forget an ambient space $\Omega$ and 
start with the following setting:  
Let $(K,d)$ be a complete metric space 
and  $\gamma = (\gamma_1,\dots , \gamma_N)$ be a system of 
proper contractionson $K$.  We assume that $K$ is self-similar, i.e.
$K=\cup_{i=1}^{N} \gamma_{i}(K)$.  We say that a system 
$\gamma = (\gamma_1,\dots , \gamma_N)$ satisfies the 
{\it open set condition in } $K$ 
if there exists an non-empty open set $V \subset K$ such that 
\[
\cup _{i=1}^N \gamma _i(V) \subset V  \text{ and } 
\gamma _i (V) \cap \gamma _j(V) = \phi 
\text{ for } i \not= j .
\]  
It is easy to see that $V$ is an open dense set of $K$. 
Moreover,  for $n \in {\mathbb N}$ and $w,v \in W_n$, if 
$w \not= v$, then $\gamma_w(V) \cap \gamma_v(V) = \phi$.  

We recall Cuntz-Pimsner algebras \cite{Pi}. See also \cite{KPW1}, 
\cite{KPW2} for the notation.
Let $A$ be a $C^*$-algebra 
and $X$ be a Hilbert right $A$-module.  We denote by $L(X)$ be 
the algebra of the adjointable bounded operators on $X$.  For 
$\xi$, $\eta \in X$, the "rank one" operator $\theta _{\xi,\eta}$
is defined by $\theta _{\xi,\eta}(\zeta) = \xi(\eta|\zeta)$
for $\zeta \in X$. The closure of the linear span of rank one 
operators is denoted by $K(X)$.   We call that 
$X$ is a Hilbert bimodule over $A$ if $X$ is a Hilbert right  $A$-
module with a homomorphism $\phi : A \rightarrow L(X)$.  We assume 
that $X$ is full and $\phi$ is injective.

   Let $F(X) = \oplus _{n=0}^{\infty} X^{\otimes n}$
be the full Fock module of $X$ with a convention 
$X^{\otimes 0} = A$. 
 For $\xi \in X$, the creation operator
$T_{\xi} \in L(F(X))$ is defined by 
$$
T_{\xi}(a) =  \xi a  \qquad \text{and } \ 
T_{\xi}(\xi _1 \otimes \dots \otimes \xi _n) = \xi \otimes 
\xi _1 \otimes \dots \otimes \xi _n .
$$
We define $i_{F(X)}: A \rightarrow L(F(X))$ by 
$$
i_{F(X)}(a)(b) = ab \qquad \text{and } \ 
i_{F(X)}(a)(\xi _1 \otimes \dots \otimes \xi _n) = \phi (a)
\xi _1 \otimes \dots \otimes \xi _n 
$$
for $a,b \in A$.  The Cuntz-Toeplitz algebra ${\mathcal T}_X$ 
is the $C^*$-algebra on $F(X)$ generated by $i_{F(X)}(a)$
with $a \in A$ and $T_{\xi}$ with $\xi \in X$.  
Let $j_K : K(X) \rightarrow {\mathcal T}_X$ be the homomorphism 
defined by $j_K(\theta _{\xi,\eta}) = T_{\xi}T_{\eta}^*$. 
We consider the ideal $I_X := \phi ^{-1}(K(X))$ of $A$. 
Let ${\mathcal J}_X$ be the ideal of ${\mathcal T}_X$ generated 
by $\{ i_{F(X)}(a) - (j_K \circ \phi)(a) ; a \in I_X\}$.  Then 
the Cuntz-Pimsner algebra ${\mathcal O}_X$ is the 
the quotient ${\mathcal T}_X/{\mathcal J}_X$ . 
Let $\pi : {\mathcal T}_X \rightarrow {\mathcal O}_X$ be the 
quotient map.  Put $S_{\xi} = \pi (T_{\xi})$ and 
$i(a) = \pi (i_{F(X)}(a))$. Let
$i_K : K(X) \rightarrow {\mathcal O}_X$ be the homomorphism 
defined by $i_K(\theta _{\xi,\eta}) = S_{\xi}S_{\eta}^*$. Then 
$\pi((j_K \circ \phi)(a)) = (i_K \circ \phi)(a)$ for $a \in I_X$.   
We note that  the Cuntz-Pimsner algebra ${\mathcal O}_X$ is 
the universal $C^*$-algebra generated by $i(a)$ with $a \in A$ and 
$S_{\xi}$ with $\xi \in X$  satisfying that 
$i(a)S_{\xi} = S_{\phi (a)\xi}$, $S_{\xi}i(a) = S_{\xi a}$, 
$S_{\xi}^*S_{\eta} = i((\xi | \eta)_A)$ for $a \in A$, 
$\xi, \eta \in X$ and $i(a) = (i_K \circ \phi)(a)$ for $a \in I_X$.
We usually identify $i(a)$ with $a$ in $A$.  We denote by 
${\mathcal O}_X^{alg}$ the $\ ^*$-algebra generated algebraically 
by $A$  and $S_{\xi}$ with $\xi \in X$. There exists an action 
$\alpha : {\mathbb R} \rightarrow Aut \ {\mathcal O}_X$
with $\alpha_t(S_{\xi}) = e^{it}S_{\xi}$, which is called the  
gauge action. Since we assume that $\phi: A \rightarrow L(X)$ is 
isometric, there is an embedding $\phi _n : L(X^{\otimes n})
 \rightarrow L(X^{\otimes n+1})$ with $\phi _n(T) = 
T \otimes id_X$ for $T \in L(X^{\otimes n})$ with the convention 
$\phi _0 = \phi : A \rightarrow L(X)$.  We denote by ${\mathcal F}_X$
the $C^*$-algebra generated by all $K(X^{\otimes n})$, $n \geq 0$ 
in the inductive limit algebra $\varinjlim L(X^{\otimes n})$. 
Let ${\mathcal F}_n$ be the $C^*$-subalgebra of ${\mathcal F}_X$ generated by 
$K(X^{\otimes k})$, $k = 0,1,\dots, n$, with the convention 
${\mathcal F}_0 = A = K(X^{\otimes 0})$.  Then  ${\mathcal F}_X = 
\varinjlim {\mathcal F}_n$.

We shall consider the union 
\[
{\mathcal G} = {\mathcal G}(\{\gamma_j:j=1,2,..,N \}) := 
\cup _{i=1}^N \{(x,y) \in K^2 ; x = \gamma _i(y)\}
\] 

of the cograph of $\gamma _i$. For example, if  
$\{\gamma_j:j=1,2,..,N \}$ 
are the continuous branches of the inverse of an expansive map 
$h: K \rightarrow K$, then ${\mathcal G}$ is exactly 
the graph of $h$. Consider a $C^*$-algebra 
$A = C(K)$ and let $X = C({\mathcal G})$.  
Then $X$ is an $A$-$A$ bimodule by 
$$
(a\cdot f \cdot b)(x,y) = a(x)f(x,y)b(y)
$$
for $a,b \in A$ and $f \in X$. We introduce an $A$-valued 
inner product $(\ |\ )_A$ on $X$ by 
$$
(f|g)_A(y) = \sum _{i=1}^N  
\overline{f(\gamma _i(y),y)}g(\gamma _i(y),y) 
$$
for $f,g \in X$ and $y \in K$. It is clear that 
the  $A$-valued inner product $(\ | \ )_A$ is well defined, 
that is, $K \ni y \mapsto (f|g)_A(y) 
\in {\mathbb C}$ is continuous. 
Put $\|f\|_2 = \|(f|f)_A\|_{\infty}^{1/2}$. 
The left multiplication of $A$ on $X$ gives 
the left action $\phi : A \rightarrow L(X)$ 
such that 
$(\phi (a)f)(x,y) = a(x)f(x,y)$ for $a \in A$ 
and $f \in X$. 

For any natural number $n$, 
we define ${\mathcal G}_n = {\mathcal G}(\{\gamma_w;w \in W_n \})$ 
and Hilbert $A$-$A$ bimodule $X_n = C({\mathcal G}_n)$ similarly. 
We also need to introduce a modified path space ${\mathcal P}_n$ 
of length $n$ by 
\[
\begin{split}
{\mathcal P}_n & = \{(\gamma _{w_1,\dots,w_n}(y),\gamma _{w_2,\dots,w_n}(y), 
\gamma _{w_3,\dots,w_n}(y), \dots, \gamma _{w_n}(y),y) \in K^{n+1} ; \\
&  w = (w_1,\dots,w_n) \in W_n , y \in K \}
\end{split}
\]
Then similarly $Y_n := C({\mathcal P}_n)$ is a $A$-$A$ bimodule  
with an $A$-valued inner product defined by 
\[
(f|g)_A(y) = \sum _{w \in W_n}   
\overline{f(\gamma _{w_1,\dots,w_n}(y),\dots, \gamma _{w_n}(y),y)}
g(\gamma _{w_1,\dots,w_n}(y),\dots, \gamma _{w_n}(y),y) 
\]
for $f,g \in Y_n$ and $y \in K$. 

If there exists a continuous function $h:K \rightarrow K$
such that each contraction $\gamma _i$ is a continuous branch 
of the inverse of $h$, 
then ${\mathcal P}_n $ can be identified with ${\mathcal G}_n $ . 
Many examples in our paper have such functions $h$.  

\begin{prop} Let $\gamma = (\gamma_1,\dots , \gamma_N)$
 be a system of 
proper contractions on a compact metric space $K$. 
Let $K$ be self-similar.
Then $X = C(\mathcal G)$  is a full Hilbert bimodule over 
$A = C(K)$ without completion. 
The left action $\phi : A \rightarrow L(X)$ is unital and 
faithful.  
Similar statements hold for $Y_n = C({\mathcal P}_n)$.
\end{prop}

\begin{proof}
For any $f \in X =  C(\mathcal G) $, we have 
$$
\| f\|_{\infty} \leq \| f \|_2 
= (\sup _y \sum _{i=1}^N |f(\gamma _i(y),y)|^2)^{1/2} 
\leq \sqrt{N} \| f\|_{\infty}
$$
Therefore two norms $\|\ \|_2$ and $\|\ \|_{\infty}$ are 
equivalent.  Since $C(\graph R)$ is complete  with  respect to 
$\| \ \|_{\infty}$, it is also complete with  respect to 
$\|\ \|_2$.

Since $(1_X|1_X)_A(y) = \sum _{i=1}^N 1_A  = N$, 
$(X|X)_A$ contains the identity $I_A$ of $A$.  Therefore 
$X$ is full. If $a \in A$ is not zero, then there exists 
$x_0 \in K$ with $a(x_0) \not = 0$. Since $K$ is self-similar, 
there exists $j$ and $y_0 \in K$ with $x_0 = \gamma _j(y_0)$. 
Choose $f \in X$ with $f(x_0, y_0) \not= 0$.  Then 
$\phi (a)f \not= 0$.  Thus $\phi$ is faithful.  
The statements for $Y_n$ are similarly proved.    
\end{proof} 

\noindent
{\bf Definition.} Let $(K,d)$ be a compact metric space 
and $\gamma = (\gamma_1,\dots , \gamma_N)$ be a system of proper contractions
on $K$.  Assume that $K$ is self-similar. 
We associate a 
$C^*$-algebra ${\mathcal O}_{\gamma}(K)$
with them as a Cuntz-Pimsner algebra ${\mathcal O}_X$ of 
the Hilbert bimodule 
$X= C({\mathcal G})$ over $A = C(K)$. 

\begin{prop}Let $(K,d)$ be a compact metric space 
and  $\gamma = (\gamma_1,\dots , \gamma_N)$ be a system of proper contractions
on $K$.  Assume that $K$ is self-similar. 
Then there exists an 
isomorphism 
$\varphi : X^{\otimes n} \rightarrow C({\mathcal P}_n)$
as a Hilbert bimodule over $A$ such that 
\begin{align*}
& (\varphi (f_1 \otimes \dots \otimes f_n))
(\gamma _{w_1,\dots,w_n}(y),\gamma _{w_2,\dots,w_n}(y), 
\gamma _{w_3,\dots,w_n}(y), \dots, \gamma _{w_n}(y),y)) \\
& = f_1(\gamma _{w_1,\dots,w_n}(y),\gamma _{w_2,\dots,w_n}(y)
    f_2(\gamma _{w_2,\dots,w_n}(y),\gamma _{w_3,\dots,w_n}(y))
\dots f_n(\gamma _{w_n}(y),y)
\end{align*}
for $f_1,\dots, f_n \in X$, $y \in K$ and 
$w = (w_1,\dots,w_n) \in W_n$.  
Moreover, let $\rho : {\mathcal P}_n \rightarrow {\mathcal G}_n$ 
be an onto continuous map such that 
\[
\rho (\gamma _{w_1,\dots,w_n}(y),\gamma _{w_2,\dots,w_n}(y), 
\gamma _{w_3,\dots,w_n}(y), \dots, \gamma _{w_n}(y),y)  
= (\gamma _{w_1,\dots,w_n}(y),y) .
\]
Then  $\rho ^* : C({\mathcal G}_n)\ni f \mapsto f \circ \rho 
\in  C({\mathcal P}_n)$ 
is an embedding as a Hilbert submodule preserving inner product. 
\label{prop:Xotimesn}
\end{prop}

\begin{proof}
It is easy to see that $\varphi$ is well-defined
and a bimodule homomorphism. We show that $\varphi$ preserves
inner product. 
Consider the case when $n = 2$ for simplicity of the notation.
\begin{align*}
& (f_1 \otimes f_2 | g_1 \otimes g_2)_A(y) 
 = (f_2|(f_1|g_1)_Ag_2)_A(y) \\
& = \sum _i
 \overline{f_2(\gamma_i(y),y)}(f_1|g_1)_A(\gamma_i(y))g_2(\gamma_i(y),y) \\
& = \sum _i
 \overline{f_2(\gamma_i(y),y)}
(\sum _j
\overline{f_1(\gamma_j \gamma_i(y),\gamma_i(y))}g_1(\gamma_j \gamma_i(y),\gamma_i(y))
g_2(\gamma_i(y),y) \\
& = \sum _{i,j}
\overline{f_1(\gamma_j \gamma_i(y),\gamma_i(y))f_2(\gamma_i(y),y)}
g_1(\gamma_j \gamma_i(y),\gamma_i(y)) g_2(\gamma_i(y),y)\\
& = \sum _{i,j}
\overline{(\varphi (f_1 \otimes f_2))(\gamma_j \gamma_i(y),\gamma_i(y),y)} 
(\varphi (g_1 \otimes g_2))(\gamma_j \gamma_i(y),\gamma_i(y),y) \\
& = (\varphi (f_1 \otimes f_2) | \varphi (g_1 \otimes g_2))(y) 
\end{align*}  

Since $\varphi$ preserves inner product, $\varphi$ is one to one. 
The non-trivial one  is to show that $\varphi$ is onto. Since 
$\varphi(1_X\otimes \dots \otimes 1_X) = 1_X$ and 
\[
\varphi (f_1 \otimes \dots \otimes f_n)
\varphi (g_1 \otimes \dots \otimes g_n)
 = \varphi (f_1g_1 \otimes \dots \otimes f_ng_n),
\]
the image of $\varphi$ is a unital 
$\ ^*$-subalgebra of $C({\mathcal P}_n)$. 
If 
\[
(\gamma _{w_1,\dots,w_n}(y), \dots, \gamma _{w_n}(y),y) \not= 
(\gamma _{u_1,\dots,u_n}(z), \dots, \gamma _{u_n}(z),z)
\]
for some $w,u \in W_n$ and $y,z \in K$, 
then there exists a certain $i$ with $1 \leq i \leq n$ such that
$\gamma _{w_i,\dots,w_n}(y) \not= \gamma _{u_i,\dots,u_n}(z)$, or 
$y \not= z$
Hence there exists $f_i \in X$ such that 
\[
f_i(\gamma _{(w_i, \dots, w_n)}(y),\gamma _{(w_{i+1}, \dots, w_n)}(y)) 
\not= f_i(\gamma _{(u_i, \dots ,u_n)}(z),
\gamma _{(u_{i+1}, \dots ,u_n)}(z)),   
\]
where for $i = n$, this means that 
$f_n(\gamma _{w_n}(y),y) \not= f_n(\gamma _{u_n}(z),z)$. 
Then 
\begin{align*}
&  \varphi (1_X \otimes \dots f_i \dots \otimes 1_X)
((\gamma _{w_1,\dots,w_n}(y),\gamma _{w_2,\dots,w_n}(y), 
\dots, \gamma _{w_n}(y),y))) \\
& \not= 
\varphi (1_X \otimes \dots f_i \dots \otimes 1_X)
(\gamma _{u_1,\dots,u_n}(z),\gamma _{u_2,\dots,u_n}(z), 
\dots, \gamma _{u_n}(z),z)).
\end{align*}
Thus the image of $\varphi$ separates the two points.  
By the Stone-Wierstrass Theorem, the image of $\varphi$ is dense in 
$C({\mathcal P}_n)$ with respect to $\| \ \|_{\infty}$ .   
Since two norms $\|\ \|_2$ and $\|\ \|_{\infty}$ are equivalent 
and $\varphi$ is isometric with respect to $\|\ \|_2$, 
$\varphi$ is onto. The rest is clear. 

\end{proof}  

{\bf Definition.}
Consider a (branched) overing map  $\pi : {\mathcal G} \rightarrow K$  defined 
by 
$\pi (x,y) = y$ for $(x,y) \in {\mathcal G}$. Define a set 
\[
B(\gamma_1,\dots,\gamma_N) := \{x \in K ; x = \gamma_i(y) = \gamma_j(y) 
\text{ for some } y \in K \text{ and } i \not= j  \}.  
\] 
Then  $B := B(\gamma_1,\dots,\gamma_N)$ is a closed set, because 
\[
B = \cup _{i \not= j} \{x \in \gamma _i(K) \cap \gamma _j(K) ; 
    \gamma _i^{-1}(x) = \gamma _j^{-1}(x) \}. 
\] 
The set $B$ is something like a branched set in 
the case of rational function  and will be described by the ideal 
$I_X := \phi ^{-1}(K(X))$ of $A$ as in \cite{KW}. 
We define a branch index $e(x,y)$ at $(x,y) \in {\mathcal G}$ by 
\[
e(x,y) := \ ^{\#} \{i \in  \{1,\dots, N\} ; \gamma_i(y) = x  \}
\]
Hence $x \in B(\gamma_1,\dots,\gamma_N)$ if and only 
if there exists $y \in K$ with $e(x,y) \geq 2$ .  For 
$x \in K$ we define 
\[
I(x) := \{i \in \{1,\dots, N\} ; \text{ there exists } y \in K 
\text{ such that } x = \gamma _i(y) \}.
\]    

\begin{lemma} 
In the above situation, if  $x \in K\setminus 
B(\gamma_1,\dots,\gamma_N)$, 
then there exists an open neighbourhood $U_x$ of $x$ satisfying the 
following: 
\begin{enumerate} 
\item $U_x \cap B = \phi$. 
\item If  $i \in I(x)$, then 
$\gamma_j(\gamma_i^{-1}(U_x)) \cap U_x = \phi $  for $j \not= i$.   
\item If $i \not\in I(x)$,  then $U_x \cap \gamma _i(K) = \phi$. 
\end{enumerate}
\label{lemma:U_x}
\end{lemma}
 
\begin{proof} Let $x \in K\setminus B$. 
Since $B$ and $\cup _{i \not\in I(x)} \gamma _i(K)$ 
are closed and $x$ is not in either of them, there exists 
an open neighbourhood $W_x$ of $x$  such that 
\[
W_x \cap (B \cup \cup _{i \not\in I(x)} \gamma _i(K)) = \phi .
\]
For $i \in I(x)$ there exists a unique $y_i \in K$ with $x = \gamma _i(y_i)$, 
since $x \not\in B$.  For $j \in \{1,\dots ,N\}$, if $j \not= i$, then 
$\gamma _j(y_i) \not= \gamma _i(y_i) = x$.  Therefore there exists an 
open neighbourhood $V_x^i$ of $x$ such that 
$\gamma_j(\gamma_i^{-1}(V_x^i)) \cap V_x^i = \phi $  for $j \not= i$.
Put $U_x := W_x \cap (\cap _{i \in I(x)} V_x^i)$.  Then 
$U_x$ is an open neighbourhood of $x$ and satisfies all the requirement.  

\end{proof}

\begin{prop}
Let $(K,d)$ be a compact metric space 
and  $\gamma = (\gamma_1,\dots , \gamma_N)$ be a system of proper contractions
on $K$.  Assume that $K$ is self-similar and 
the system $\gamma = (\gamma_1,\dots , \gamma_N)$ satisfies the open 
set condition in $K$.
Then  
\[
I_X = \{a \in A=C(K) ; a \text{ vanishes on } B(\gamma_1,\dots,\gamma_N) \}.
\]
\label{prop:critical}
\end{prop}
\begin{proof} 
Let $B = B(\gamma_1,\dots,\gamma_N)$. 
Firstly, let us take  $a \in A$ with a compact support 
$S = \supp (a) \subset K \setminus B$. For any $x \in S$, 
choose an open neighbourhood $U_x$ of $x$ as in Lemma \ref{lemma:U_x}. 
Since $S$ is compact, there exists a finite 
subset $\{x_1,\dots , x_m\}$ such that 
$S \subset \cup _{i=1}^m U_{x_i} \subset K \setminus B$. 
By considering a partition of unity for an open covering 
$K = S^c \cup \cup _{i=1}^m U_{x_i}$, 
we can choose  a finite family $(f_i)_i$ in $C(K)$ such that 
$0 \leq f_i \leq 1$, $\supp(f_i) \subset U_{x_i}$ for 
$i=1,\dots ,m$ and $\sum _{i=1}^{m} f_i(x) = 1$ 
for $x \in S$.  Define 
$\xi _i, \eta _i \in C({\mathcal G})$ by 
$\xi _i(x,y) = a(x)\sqrt{f_i(x)}$ and 
$\eta _i(x,y) = \sqrt{f_i(x)}$.  Consider 
$T:= \sum _{i = 1}^k \theta _{\xi _i, \eta _i} \in K(X)$.
We shall show that $T = \phi (a)$.      
For any $\zeta \in C({\mathcal G})$, we have 
$(\phi (a)\zeta )(x,y) = a(x)\zeta (x,y)$ and 
\begin{align*}
& (T\zeta )(x,y) = \sum _i \xi _i(x,y) 
\sum _j
\overline{\eta _i (\gamma_j(y),y)}\zeta (\gamma_j(y),y) \\
& = \sum _i a
(x)\sqrt{f_i(x)} 
\sum _j \sqrt{f_i(\gamma_j(y))}\zeta (\gamma_j(y),y).
\end{align*}
In the case when $a(x) = 0$, we have 
\[
(T\zeta )(x,y) = 0 = (\phi (a)\zeta )(x,y).
\]
In the case when $a(x) \not= 0$, we have $x \in \supp (a) = S 
\subset \cup _{i=1}^m U_{x_i}$. Hence 
$x \in U_{x_i}$ for some $i$. Take any $y \in K$ with $(x,y) \in 
{\mathcal G}$.  Since $x \not\in B$, there
exists a unique  
$k \in \{1,\dots, N\}$ with $x = \gamma_k(y)$.  
Then for any $j \not=k$ 
$f_i(\gamma_j(y)) = 0$, because 
$\gamma_j(y) \in \gamma_j(\gamma_k^{-1}(U_{x_i}) \subset U_{x_i}^c$. 
Therefore we have 
\begin{align*}
(T\zeta )(x,y) & = \sum _i a
(x)\sqrt{f_i(x)} 
(\sum _j \sqrt{f_i(\gamma_j(y))}\zeta (\gamma_j(y),y))\\
& = \sum _i a
(x)\sqrt{f_i(x)} \sqrt{f_i(\gamma_k(y))}\zeta (\gamma_k(y),y) \\
& = \sum _i a(x)f_i(x)\zeta (x,y) 
= a(x)\zeta (x,y) = (\phi (a)\zeta )(x,y).
\end{align*}
Thus $\phi (a) = T \in K(Y)$.  Now for a general $a
\in A$ which vanishes on $B$, there exists a sequence 
$(a_n)_n$ in $A$ with 
compact supports $\supp (a_n) \subset  
K \setminus B$ such that 
$\|a - a_n\|_{\infty} \rightarrow 0$.  
Hence $\phi (a) \in K(X)$, i.e., $a \in I_X$.

Conversely let $a \in A$ and $a(c) \not= 0$ for some 
$c \in B$.  We may assume that $a(c) = 1$. 
Then $c = \gamma_k(d) = \gamma_r(d)$  
for some $d \in K$ with $k \not= r \in \{1,\dots, N\}$. 
Thus the branch index $e(c,d) \geq 2$.  
We need to show that $\phi(a) \not\in K(X)$.  
On the contrary suppose that $\phi(a) \in K(X)$.  
Then for $\varepsilon = \frac{1}{5\sqrt{N}}$, there exists a finite subset 
$\{\xi _i, \eta _i \in X ; i = 1,\dots ,M\}$ 
such that $\| \phi(a) - 
\sum _{i =1}^M \theta _{\xi _i, \eta _i} \| < \varepsilon$.  
Since the system satisfies the open set condition in $K$, 
there exists an open dense set $V \subset K$ such that 
$\cup _{i=1}^N \gamma (V) \subset V$ and  
$\gamma _i (V) \cap \gamma _j(V) = \phi $ 
for $i \not= j$ .  Thus $\gamma_i(V)$ is dense in 
$\gamma_i(K)$ and 
${\mathcal G}_V 
:= \cup _{i=1}^{N}\{(\gamma_i(y),y) \in {\mathcal G} ; y \in V\}$ 
is dense in ${\mathcal G}$. 
We claim that for any open neighbourhood $U_{(c,d)}$ of $(c,d)$ in 
${\mathcal G}$, there exists $(x,y) \in U_{(c,d)}$ with $e(x,y) = 1$.  
On the contrary suppose that there were an 
open neighbourhood $U_{(c,d)}$ of $(c,d)$ in ${\mathcal G}$ such that 
for any  $(x,y) \in U_{(c,d)}$ we have  $e(x,y) \geq 2$. Then 
there exists $(x,y) \in {\mathcal G}_V \cap U_{(c,d)}$ with $e(x,y) \geq 2$.
Thus $y \in V$ and there exist $i$ and $j$ such that  $i \not= j$ and 
$x = \gamma_i(y) = \gamma_j(y)$. Then 
$x \in \gamma _i(V) \cap \gamma _j(V)$. 
This is a contradiction and the claim is shown. Therefore there exists 
a sequence $(x_n,y_n)_n$ in ${\mathcal G}$ such that $e(x_n,y_n) = 1$ and 
$(x_n,y_n)_n$ converges to $(c,d)$.  Since ${\mathcal G}$ is a finite union 
of $\{(\gamma_i(y),y) ; y \in K\}$, $i = 1,\dots, N$, we may assume that 
there exists a certain $i_0$ such that 
$\{(x_n,y_n) ; n \in {\mathbb N}\} 
\subset \{(\gamma_{i_0}(y),y) ; y \in K\}$, 
by taking a subsequence if necessary.  Since 
$e(x_n,y_n) = 1$, as in the proof of  Lemma \ref{lemma:U_x}, 
there exists an open neighbourhood $U_n = U_{x_n}$ of $x_n$ 
such that
$\gamma_j(\gamma_{i_0}^{-1}(U_n)) \cap U_n= \phi$ 
for $j \not= i_0$.    
We choose  $\zeta _n \in X$ such that 
$\supp \zeta _n  \subset 
\{(x,y) \in {\mathcal G} ; x \in U_n \text{ and } x = \gamma _{i_0}(y) \}$, 
$\zeta _n(x_n,y_n) = 1$ and 
$0 \leq \zeta _n \leq 1$.  Then 
$\|\zeta _n \|_2 \leq \sqrt{N}$. 
If $j \not= i_0$, then $\gamma_j(y_n) \not\in U_n$  
and $\zeta _n(\gamma_j(y_n),y_n) = 0$.  
If $j = i_0$, then $\zeta _n(\gamma_{i_0}(y_n),y_n) 
= \zeta _n(x_n,y_n) = 1$.  
Hence 
\begin{align*}
& |a(x_n) - \sum _{i=1}^{M} \xi _i(x_n,y_n)
\overline{\eta _i(x_n,y_n)}| \\
& = |a(x_n) - 
\sum _{i=1}^{M} \xi _i(x_n,y_n) 
\sum _{j=1}^N 
\overline{\eta _i (\gamma_j(y_n),y_n)}\zeta _n(\gamma_j(y_n),y_n)|\\
& = |((\phi(a) - \sum _{i =1}^M \theta _{\xi _i, \eta _i})
  \zeta _n)(x_n,y_n)| \\
& \leq \|(\phi(a) - \sum _{i =1}^M \theta _{\xi _i, \eta _i})
  \zeta _n\|_2 
\leq \|(\phi(a) - \sum _{i =1}^M \theta _{\xi _i, \eta _i}\|
  \|\zeta _n\|_2 
\leq \varepsilon \sqrt{N}.
\end{align*}  
Since $(x_n,y_n) \rightarrow (c,d)$ as $n \rightarrow \infty$, 
we have 
\[
|a(c) - \sum _{i=1}^{M} \xi _i(c,d)
\overline{\eta _i(c,d)}| 
\leq \varepsilon \sqrt{N}.
\]
On the other hand, consider  $\zeta \in X$ satisfying 
$\zeta (c,d) = 1$,  $0 \leq \zeta \leq 1$ and 
$\zeta (\gamma_j(d),d) = 0$ 
for $j$ with $\gamma_j(d) \not= c$. 
Then 
\begin{align*}
& |a(c) - \sum _{i=1}^{M} \xi _i(c,d)
e(c,d)\overline{\eta _i(c,d)}| \\
& = |a(c) - \sum _{i=1}^{M} \xi _i(c,d)
\sum_{j=1}^N \overline{\eta _i
(\gamma_j(d),d)}\zeta (\gamma_j(d),d)| \\
& \leq \|(\phi(a) - \sum _{i =1}^M \theta _{\xi _i, \eta _i})
  \zeta \|_2 
\leq \varepsilon \sqrt{N}.
\end{align*}  

Since $e(c,d) \geq 2$ and $a(c) = 1$,  we have 
\begin{align*}
\frac{1}{2} &
\leq |a(c) - \frac{1}{e(c,d)}a(c)| \\
& \leq |a(c) -  \sum _{i=1}^{M} \xi _i(c,d) \overline{\eta _i(c,d)}|  
- |\sum _{i=1}^{M} \xi _i(c,d) 
\overline{\eta _i(c,d)}-\frac{1}{e(c,d)}a(c)| \\
& \leq \varepsilon \sqrt{N} + \frac{1}{e(c,d)}\varepsilon \sqrt{N} 
\leq 2\varepsilon \sqrt{N} = \frac{2}{5}
\end{align*}

This is a contradiction. 
Therefore  $\phi(a) \not\in K(X)$
\end{proof}

\begin{cor}
$\ ^{\#}B(\gamma_1,\dots,\gamma_N) = \dim (A/I_X)$ 
\end{cor}

\begin{cor} The closed set $B(\gamma_1,\dots,\gamma_N) = \phi$ 
if and only if 
$\phi (A)$ is contained in $K(X)$ if and only if
$X$ is finitely generated projective right $A$ module.
\end{cor}

\section{Simplicity and pure infinteness}
Let $(K,d)$ be a compact metric space 
and $\gamma = (\gamma_1,\dots , \gamma_N)$ be a system of proper contractions
on $K$.  Assume that $K$ is self-similar. 
Let $A = C(K)$ and $X = C(\mathcal G)$. 
Define an endomorphism $\beta _i : A \rightarrow A$ by 
\[
(\beta _i (a))(y) = a(\gamma _i (y))
\]
for $a \in A$, $y \in K$.  We also define a unital completely 
positive map $E_{\gamma}: A \rightarrow A$ by
\[
(E_{\gamma}(a))(y) := \frac{1}{N} 
\sum _{i=1}^N a(\gamma _i (y))
\]
for $a \in A$, $y \in K$, that is, 
$E_{\gamma} =  \frac{1}{N} \sum _{i=1}^N \beta _i$.  
For a constant function 
$\xi _0 \in X$ with 
\[
\xi _0 (x,y) := \frac{1}{\sqrt{N}}
\]
we have     
\[
E_{\gamma}(a) = (\xi _0|\phi (a)\xi _0)_A  
\text{ and } \ E_{\gamma}(I) = (\xi _0|\xi _0)_A = I. 
\]
We introduce an operator $D := S_{\xi _0} \in 
{\mathcal O}_{\gamma}(K)$.

\begin{lemma} 
In the above situation, for  $a\in A$, we have the following:
\[
D^*aD = E_R(a) \text{ and in particular } D^*D = I. 
\]

\end{lemma}
 
\begin{proof}
\[
D^*aD = S_{\xi _0}^*aS_{\xi _0} = (\xi _0|\phi (a)\xi _0)_A 
= E_R(a).
\]

\end{proof}

{\bf Definition}. Let $(K,d)$ be a complete metric space 
and  $\gamma = (\gamma_1,\dots , \gamma_N)$ be a system of proper contractions
on $K$. Then $a \in A = C(K)$ is said to be 
$\gamma$-invariant if 
\[
a(\gamma _i (y)) = a(\gamma _j (y))   \ \ \ \  
\text{ for any  } y \in K  \text { and } i,j = 1, \dots , N
\]
Suppose that  $K$ is a self-similar set and 
$a \in A = C(K)$ is 
$(\gamma_w)_{w \in W_n}$-invariant, 
then $a$ is $(\gamma_w)_{w \in W_{n-1}}$-invariant. 
In fact, for any $y \in K$  
there exists $z \in K$ and $i$ such that 
$y = \gamma _i(z)$, since $K$ is self-similar.  
Then for any  $w, v \in W_{n-1}$, we have 
\[
 a(\gamma _w (y)) = a(\gamma _{wi} (z)) 
= a(\gamma _{vi} (z)) = a(\gamma _v (y))
\]
If $a$ is  $(\gamma_w)_{w \in W_n}$-invariant, then 
for any $k = 1, \dots, n$ we may write 
\[
\beta ^k(a)(y) := a(\gamma _{w_1} \dots \gamma _{w_k}(y)), 
\ \text{ for any } w \in W_k .
\]
Since  $\beta ^k(a)(y)$ does not depend on the choice of 
$w \in W_k$, $\beta ^k(a)(y)$ is well defined. We may 
write  that $\beta (\beta ^{k-1}(a))(y) = \beta ^k(a)(y)$.

\begin{lemma} 
In the same situation, if $a \in A$
is  $(\gamma_w)_{w \in W_n}$-invariant,
then for any 
$f_1,\dots , f_n \in X$, 
we have the following: 
$$
aS_{f_1}\dots S_{f_n} = S_{f_1}\dots S_{f_n}\beta ^n(a) .
$$
\label{lemma:a-f}
\end{lemma}

\begin{proof} If $a \in A$
is  $(\gamma_w)_{w \in W_n}$-invariant,
 then $\beta (a)$ is $(\gamma_w)_{w \in W_{n-1}}$-invariant.  
Therefore it is enough to show that 
$aS_f = S_f\beta (a)$ for $f \in X$. 
We have  $aS_f = S_{\phi (a)f}$ 
and $S_f\beta (a) = S_{f\beta (a)}$.  
Since 
\begin{align*}
f\beta (a)(\gamma _j(y),y) 
& = f(\gamma _j(y),y)(\beta (a))(y)
  = f(\gamma _j(y),y)a(\gamma _i(y)) \\
& = a(\gamma _j(y))f(\gamma _j(y),y) 
  = (\phi (a)f)(\gamma _j(y),y)  
\end{align*}
 
\noindent
we have $aS_f = S_f\beta (a)$. 
\end{proof}

\par
\begin{lemma} Let $(K,d)$ be a compact metric space 
and $\gamma = (\gamma_1,\dots , \gamma_N)$ be a system of proper contractions
on $K$.  Assume that $K$ is self-similar.
For any non-zero positive element $a \in A$ and for any  
$\varepsilon > 0$ there exist $n \in \mathbb{N}$ and 
$f \in X^{\otimes n}$  with $(f|f)_A = I$ such that  
$$
 \|a\| -\varepsilon \le S_f^* a S_f \le \|a\|
$$
\label{lemma:epsilon}
\end{lemma}
\begin{proof} Let $x_0$ be a point in $K$ with 
$|a(x_0)| = \| a \|$.  For any  
$\varepsilon > 0$ there exist an open neighbourhood $U_0$ of $x_0$
in $K$ such that for any $x \in U_0$  we have 
$\|a\| -\varepsilon  \le a(x) \le \|a\| $. 
Choose  anothter open neighbourhood $U_1$ of $x_0$ in 
$K$ and a compact subset $K_1 \subset K$ satisfying 
$U_1 \subset K_1 \subset U_0$. Then there exists 
$n \in {\mathbb N}$ and $v \in W_n$ such that 
$\gamma _v(K) \subset U_1$. 
We identify $X^{\otimes n}$ with $C({\mathcal P}_n) 
\supset \rho ^*(C({\mathcal G}_n)) $
as in Proposition \ref{prop:Xotimesn}. 
Define closed subsets $F_1$ and $F_2$ of $K \times K$ by  

\begin{align*}
 F_1 & = \{(x,y)\in K \times K ; x = \gamma_w (y), 
x \in K_1 \text{ for some } w \in W_n \} \\
 F_2 & = \{(x,y) \in K \times K ; x = \gamma_w (y),
x \in U_0^{c} \text{ for some }w \in W_n \}
\end{align*}

Since $F_1 \cap F_2 = \phi $, there exists  
$g \in C({\mathcal G}_n)$ such that 
$0 \le g(x,y) \le 1$  and 
$$
 g(x,y) =\left\{\begin{array}{cc}
  1,  & (x,y)\in F_1  \\
  0,  & (x,y)\in F_2
       \end{array}\right.
$$
\par
Since 
$\gamma _v(K) \subset U_1$, for any $y \in K$ there exists 
$x_1 \in U_1$ such that $x_1 = \gamma _v(y) \in U_1 \subset K_1$, 
so that $(x_1,y) \in F_1$. Therefore  
\[
 (g|g)_A(y) = \sum_{w \in W_n} |g(\gamma_w(y),y)|^2 
             \ge |g(x_1,y)|^2  \ge 1 .
\]

Let $b:=(g|g)_A$.  Then $b(y) = (g|g)_A(y) \ge 1$.  Thus 
$b \in A$ is positive and invertible.
We put $f :=\rho ^*(g b^{-1/2}) = \rho ^*(g)b^{-1/2}\in X^{\otimes n}$.  Then  
$$
 (f|f)_A   = (g b^{-1/2}| g b^{-1/2})_A 
          = b^{-1/2} (g|g)_A b^{-1/2} 
           = I.
$$

\par
For any $y \in K$ and any $w=(w_1,\dots,w_n) \in W_n$,  
let  $x = \gamma_w(y)$. 
If $x \in U_0$, 
then $\|a\|-\varepsilon \le a(x)$, and  
if $x \in U_0^c$, then  
\[
f(\gamma _{w_1,\dots,w_n}(y),\dots, \gamma _{w_n}(y),y)
=g(x,y) b^{-1/2}(y) = 0,
\] 
because $(x,y) \in F_2$. 
Therefore   
\begin{align*}
 \|a \|-\varepsilon & = (\|a\|-\varepsilon) (f|f)_A(y) \\
                    & = (\|a\|-\varepsilon) \sum_{w \in W_n} 
              |f(\gamma _{w_1,\dots,w_n}(y),\dots, \gamma _{w_n}(y),y)|^2 \\
               & \leq \sum_{w \in W_n}
    a(\gamma_w(y))|f(\gamma _{w_1,\dots,w_n}(y),\dots, \gamma _{w_n}(y),y)|^2 \\
                    & = (f|af)_A(y) = S_f^* a S_f (y). 
\end{align*}

We also have that 
\[
S_f^* a S_f = (f|af)_A \leq \| a \| (f|f)_A = \| a \| .
\] 

\end{proof}

\begin{lemma}
Let $(K,d)$ be a compact metric space 
and $\gamma = (\gamma_1,\dots , \gamma_N)$ be a system of proper contractions
on $K$.  Assume that $K$ is self-similar.
For any non-zero positive element $a \in A$ and for any  
$\varepsilon > 0$ with $0 < \varepsilon < \|a\|$, 
there exist $n \in \mathbb{N}$ and $u \in X^{\otimes n}$ 
such that  
$$
 \| u \| _2 \le (\|a\| - \varepsilon)^{-1/2} \qquad \text{and} 
\quad  S_u^* a S_u = I
\label{lemma:u-epsilon}
$$
\end{lemma}
\begin{proof} 
For any $a \in A$ and $\varepsilon > 0$ as above, we choose 
$f \in X^{\otimes n}$ as in Lemma \ref{lemma:epsilon}. 
Put $c=S_f^*aS_f$.  
Since $0 < \|a\|-\varepsilon \le c \le \|a\|$, $c$ is positive 
and invertible.  Let $u:=fc^{-1/2}$. Then 
$$
 S_u^* a S_u = (u|au)_A = (fc^{-1/2}| a fc^{-1/2})_A 
             = c^{-1/2}(f|af)_A c^{-1/2} = I.
$$
Since $\|a\|-\varepsilon \le c$, 
we have $c^{-1/2}\le (\|a\|-\varepsilon)^{-1/2}$.¡¡
Hence 
$$
 \|u\| _2 = \|f c^{-1/2}\| _2 \le \|c^{-1/2}\| _2 
\leq (\|a\| - \varepsilon)^{-1/2}.
$$
\end{proof}

We need the following easy fact:
Let $F$ be a closed subset of a topological space $Z$.  Let 
$a: F \rightarrow \mathbb{C}$ be continuous.  If $a(x) = 0$ for 
$x$ in the boundary of $F$, Then $a$ can be extended to a 
continuous function on $Z$ by putting $a(x) = 0$ for 
$x \not\in F$.   

\begin{lemma}
Let $(K,d)$ be a compact metric space 
and  $\gamma = (\gamma_1,\dots , \gamma_N)$ be a system of proper contractions
on $K$.  Assume that $K$ is self-similar and 
the system $\gamma = (\gamma_1,\dots , \gamma_N)$ satisfies the open 
set condition in $K$.  
For any $n \in \mathbb{N}$,  
any $T \in L(X^{\otimes n})$ and 
any  $\varepsilon > 0$,  there exists a positive element $a \in A$ 
such that $a$ is $\{\gamma_w ; w \in W_n \}$-invariant, 
\[
 \| \phi(a)T\|^2 \ge \|T\|^2 -\varepsilon
\]
and $\beta^p(a) \beta^q(a) = 0$ for 
$p,q=1, \cdots ,n$ with $p \not= q$.
\label{lemma:T-epsilon}
\end{lemma}

\begin{proof}
For any $n \in {\mathbb N}$,  any $T \in L(X^{\otimes n})$ and 
any  $\varepsilon > 0$, there exists $f \in X^{\otimes n}$
such that $ \|f\|_2=1$ and 
$\|T\|^2 \ge \|Tf\|_2^2 > \|T\|^2 -\varepsilon $. 
We still identify $X^{\otimes n}$ with $C({\mathcal P}_n)$.  
Then there exists $y_0 \in K$ such that  
\[
 \|Tf\|_2^2 =  \sum_{w \in W_n}
 |(Tf)(\gamma _{w_1,\dots,w_n}(y_0),\dots, \gamma _{w_n}(y_0),y_0)|^2
  > \|T\|^2 -\varepsilon.
\]
Since $y \mapsto (Tf|Tf)_A(y)$ is continuous and 
\[
 \|Tf\|_2^2 = \sup_{y \in K} \sum_{w \in W_n}
|(Tf)(\gamma _{w_1,\dots,w_n}(y),\dots, \gamma _{w_n}(y),y)|^2,
\] 
there exists an open neighbourhood $U_0$ of $y_0$ such that 
for any $y \in U_0$
\[
   \sum_{w \in W_n}
|(Tf)(\gamma _{w_1,\dots,w_n}(y),\dots, \gamma _{w_n}(y),y)|^2
 > \|T\|^2 -\varepsilon .
\]
Since $\gamma = (\gamma_1,\dots , \gamma_N)$ satisfies the open 
set condition in $K$,  
there exists an open dense $V \subset K$ such that 
\[
\cup _{i=1}^N \gamma _i (V) \subset V  \text{ and } 
\gamma _i (V) \cap \gamma _j(V) = \phi 
\text{ for } i \not= j .
\]  
Then there exist $y_1 \in V \cap U_0$ and an 
open neighbourhood $U_1$ of $y_1$ with 
$U_1 \subset V \cap U_0$.    
Since the contractions are proper and $K$ is self-similar, 
there exist $r \in {\mathbb N}$ and 
$(j_1,\dots ,j_r) \in W_r$ such that 
\[
\gamma_{j_1} \gamma_{j_2} \dots \gamma_{j_r}(V) 
\subset U_1 \subset V \cap U_0 .
\]  
Put $j_{r+1} = 2$  and $j_{r+2} = j_{r+3} = \dots = j_{r+n} = 1$. 
Then 
\[
\phi \not= \gamma_{j_1} \gamma_{j_2} \dots \gamma_{j_{r+n}}(V) 
\subset \gamma_{j_1} \gamma_{j_2} \dots \gamma_{j_r}(V) 
\subset U_1 \subset V \cap U_0 .
\]  
There exist $y_2 \in K$, an open neighbourhood $U_2$ of $y_2$ 
and a compact set $L$  such that 
\[
y_2 \in U_2 \subset L \subset 
\gamma_{j_1} \gamma_{j_2} \dots \gamma_{j_{r+n}}(V) 
\subset U_1 \subset V \cap U_0 .
\] 
Choose a positive function $b \in A$ such that 
$0 \leq b \leq 1$,  $b(y_2) = 1$ and $b|_{U_2^c} = 0$. 
Thus $\{x \in K ; b(x) \not= 0 \} \subset U_2$. 
For $w \in W_n$, we have 
\[
\gamma_w(y_2) \in \gamma_w(U_2)
\subset \gamma_w(L) \subset \gamma_w(V). 
\]
Moreover for  $w,v \in W_n$, by open set condition, 
\[
\gamma_w(L) \cap \gamma_v(L) = \phi   \text{ if }
 w \not= v .
\]
Now we denine a positive function $a$ on $K$ by 
$$
 a(x) =\left\{\begin{array}{cc}
  b(\gamma_w^{-1}(x)),  & \text{ if } 
                              x \in \gamma_w(L), \ w \in W_n  \\
  0,  & \text{ if  otherwise }
       \end{array}\right.
$$ 
Since $L' := \cup _{w \in W_n} \gamma_w(L)$ is compact, 
$U' := \cup _{w \in W_n} \gamma_w(U_2)$ is open and 
$\{ x \in K ; a(x) \not= 0 \} \subset U' \subset L'$, 
$a$ is continuous on $L'$ and $a(x) = 0$ for 
$x$ in the boundary of $L'$.  Therefore $a$ is continuous 
on $K$, i.e. $a \in A = C(K)$. By the construction, 
$a$ is $(\gamma_w)_{w \in W_n}$-invariant. 

For a natural number $p \leq n$ and $(i_1,\dots,i_p) \in  W_p$, 
we have 
\[
\supp (\beta_{i_p} \beta_{i_{p-1}}\dots \beta_{i_1}(a)) 
\subset \cup _{(i_{p+1},\dots,i_n) \in W_{n-p}} 
\gamma_{i_{p+1}} \dots \gamma_{i_{n}} (\supp b) .
\]
In fact, if $a(\gamma_{i_1}\dots \gamma_{i_p}(z)) \not= 0$, then 
there exists $(i_{p+1},\dots,i_n) \in W_{n-p}$ and $y \in L$ 
satisfying $z = \gamma_{i_{p+1}} \dots \gamma_{i_n}(y)$ by the 
definitin of $a$. Moreover 
\[
a(\gamma_{i_1} \dots \gamma_{i_p}(z)) = 
b(\gamma_{(i_1,\dots,i_n)}^{-1}(\gamma_{i_1}\dots \gamma_{i_p})(z))
= b(\gamma_{(i_{p+1},\dots,i_n)}^{-1}(z)) \not= 0 .
\]
Hence $z \in \gamma_{i_{p+1}} \dots \gamma_{i_n} (\supp b)$.  

Since $\supp b \subset L \subset 
\gamma_{j_1} \gamma_{j_2} \dots \gamma_{j_{r+n}}(V)$, 
\[
\supp (\beta_{i_p} \beta_{i_{p-1}}\dots \beta_{i_1}(a)) 
\subset \bigcup_{(i_{p+1},\dots,i_n) \in W_{n-p}} 
\gamma_{i_{p+1}} \dots \gamma_{i_n} 
\gamma_{j_1} \gamma_{j_2} \dots \gamma_{j_{r+n}}(V).
\]
For $1 \leq p \lneqq q \leq n$, 
\[
\supp (\beta ^p(a)) 
\subset \bigcup_{(i_{p+1},\dots,i_n) \in W_{n-p}} 
\gamma_{i_{p+1}} \dots \gamma_{i_n} 
\gamma_{j_1} \gamma_{j_2} \dots \gamma_{j_{r+n}}(V).
\]
and 
\[
\supp (\beta ^q(a)) 
\subset \bigcup _{(i_{q+1},\dots,i_n) \in W_{n-q}} 
\gamma_{i_{q+1}} \dots \gamma_{i_n} 
\gamma_{j_1} \gamma_{j_2} \dots \gamma_{j_{r+n}}(V).
\]
Since $(n-p)+(r+1)$-th subsuffixes are different as   
$j_{r+1} = 2 \not= 1 = j_{r+1+(q-p)}$, 
we have  
$\supp (\beta ^p(a)) \cap \supp (\beta ^q(a)) = \phi$. 
Thus $\beta^p(a) \beta^q(a) = 0$

Furthermore, we have 
\begin{align*}
\| \phi(a)Tf\| _2^2  
 & =  \sup_{y \in K} \sum_{w \in W_n}
 |(a(\gamma_w(y))(Tf)(\gamma _{w_1,\dots,w_n}(y),\dots, \gamma _{w_n}(y),y)|^2  \\
 & =  \sup_{y \in L} \sum_{w \in W_n}
 |(b(y)(Tf)(\gamma _{w_1,\dots,w_n}(y),\dots, \gamma _{w_n}(y),y)|^2  \\
 & \ge  \sum_{w \in W_n}
 |(Tf)(\gamma _{w_1,\dots,w_n}(y_2),\dots, \gamma _{w_n}(y_2),y_2)b(y_2)|^2 \\
 & =    \sum_{w \in W_n}
 |(Tf)(\gamma _{w_1,\dots,w_n}(y_2),\dots, \gamma _{w_n}(y_2),y_2)|^2 \\
 & > \|T\|^2 -\varepsilon. 
\end{align*} 
because $y_2 \in L \cap U_2 \subset U_0$.  
Therefore we have $\| \phi(a)T\|^2 \ge \|T\|^2 -\varepsilon$.

\end{proof}

Let ${\mathcal F}_n$ be the $C^*$-subalgebra of ${\mathcal F}_X$ 
generated by $K(X^{\otimes k})$, $k = 0,1,\dots, n$ 
and $B_n$ be the $C^*$-subalgebra of ${\mathcal O}_X$ 
generated by 
\[
\bigcup_{k=1}^n \{S_{x_1} \dots S_{x_k}S_{y_k}^* \dots S_{y_1}^* : 
x_1, \dots x_k,  y_1, \dots y_k \in X \} \cup A.
\]
 In the following Lemma \ref{lemma:I-free} we shall use 
an isomorphism 
$\varphi : {\mathcal F}_n \rightarrow B_n$ 
as in Pimsner \cite{Pi}  and  Fowler-Muhly-Raeburn \cite{FMR}
such that 
$$
\varphi (\theta_{x_1 \otimes \dots \otimes x_k, 
         y_1 \otimes \dots \otimes y_k}) 
    =  S_{x_1} \dots S_{x_k}S_{y_k}^* \dots S_{y_1}^*.
$$
To simplify notation, we put $S_x = S_{x_1} \dots S_{x_k}$ 
for $x = x_1 \otimes \dots \otimes x_k \in X^{\otimes k}$ .

\begin{lemma}
In the above situation,  
let $b = c^*c$ for some $c \in {\mathcal O}_X^{alg}$.  
We decompose  $b = \sum _j b_j$ with 
$\gamma _t(b_j) = e^{ijt}b_j$.
For any  $\varepsilon >0 $
there exists $P \in A$ with $0\le P \le I$ satisfying the 
following:  
\begin{enumerate}
 \item $Pb_jP = 0$ \qquad $(j\ne 0)$
 \item $\|Pb_0P\| \ge \|b_0\| -\varepsilon $
\end{enumerate} 
\label{lemma:I-free}
\end{lemma} 

\begin{proof} For $x \in X^{\otimes n}$, we define $\length (x) = n$ 
with the convention $\length (a) = 0$ for $a \in A$.  
We write $c$ as a finite sum $c = a + \sum _i S_{x_i}S_{y_i}^*$.
Put $n = 2 \max \{\length (x_i), \length (y_i) ; i\}$.  

\par\noindent
For $j > 0$, each $b_j$ is a finite sum of terms in the form such that 
$$
S_x S_y^* \qquad x \in X^{\otimes (k+j)}, \qquad y \in X^{\otimes k} 
\qquad 0 \le k+j \le n 
$$
In the case when $j<0$, 
$b_j$ is a finite sum of terms in the form such that 
$$
 S_x S_y^* \qquad x \in X^{\otimes k}, \qquad y \in X^{\otimes (k+|j|)} 
\qquad 0 \le k+|j| \le n 
$$

We shall identify $b_0$ with an element in 
$A_{n/2} \subset A_n \subset L(X^{\otimes n})$.  
Apply Lemma \ref{lemma:T-epsilon} for   and 
$T = (b_0)^{1/2}$. Then  there exists
a positive element $a \in A$ such 
that $a$ is $\{\gamma_w ; w \in W_n \}$-invariant,  
$ \| \phi(a)T\|^2 \ge \|T\|^2 -\varepsilon$, 
\[
 \| \phi(a)T\|^2 \ge \|T\|^2 -\varepsilon
\]
and $\beta^p(a) \beta^q(a) = 0$ for 
$p,q=1, \cdots ,n$ with $p \not= q$.
Define a positive operator $P = a \in A$.  Then 
$$
  \| P b_0 P\| = \| Pb_0^{1/2} \|^2 
                \ge \|b_0^{1/2}\|^2 -\varepsilon 
                = \| b_0 \| -\varepsilon
$$
For $j > 0$, we have   
\begin{align*}
 PS_xS_y^* P & =  aS_xS_y^* a 
             =  S_x \beta^{k+j}(a) \beta^{k}(a) S_y^* = 0
\end{align*}
For $j<0$, we also have that $PS_xS_y^* P = 0$. 
Hence $Pb_jP = 0$ for $j \not= 0$.  

\end{proof}

\begin{thm}
Let $(K,d)$ be a compact metric space 
and  $\gamma = (\gamma_1,\dots , \gamma_N)$ be a system of proper contractions
on $K$.  Assume that $K$ is self-similar and 
the system $\gamma = (\gamma_1,\dots , \gamma_N)$ satisfies the open 
set condition in $K$.  
Then the associated $C^*$-algebra 
${\mathcal O}_{\gamma}(K)$ 
is simple and purely infinite.
  
\end{thm}

\begin{proof}
Let $w \in {\mathcal O}_X = 
{\mathcal O}_{\gamma}(K)$ 
be any non-zero positive element.    
We shall show that there exist $z_1$, 
$z_2 \in {\mathcal O}_{\gamma}(K)$  
such that $z_1^*w z_2 =I$.  
We may assume that $\|w\|=1$.
Let $E : {\mathcal O}_{\gamma}(K)
 \rightarrow {\mathcal O}_{\gamma}(K)^{\alpha}$
be the canonical conditional expectation onto the fixed point 
algebra by the gauge action $\alpha$. 
Since $E$ is faithful, $E(w) \not= 0$.  
Choose  $\varepsilon$ such that 
\[
0 < \varepsilon < \frac{\|E(w)\|}{4} \  \text{ and } \ 
\varepsilon \|E(w) -3\varepsilon \|^{-1} \leq 1 .  
\]

There exists an element $c \in {\mathcal O}_X^{alg}$  
such that $ \|w - c^* c\| < \varepsilon$ and  
$\|c\| \le 1$.  Let $b = c^*c$. Then $b$ is decomposed  
as a finite sum $b = \sum_j b_j$ with 
$\alpha_t(b_j) =e^{ijt}b_j$.  
Since  $\|b\| \le 1$,  $\|b_0\| = \|E(b)\| \le 1$. 
By Lemma \ref{lemma:I-free}, there exists $P \in A$ with 
$0 \le P \le I$ satisfying $Pb_jP = 0$ \qquad  $(j\ne 0)$
and $\|Pb_0P\| \ge \|b_0\| -\varepsilon $. 
Then we have 

\begin{align*}
\| Pb_0 P \|  & \ge \|b_0 \| -\varepsilon 
              = \|E(b)\| -\varepsilon \\
              &  \ge \|E(w)\| -\|E(w) - E(b)\| -\varepsilon 
              \ge \|E(w)\| -2 \varepsilon 
\end{align*}
For $T := Pb_0P \in L(X^{\otimes m})$,  
there exists $f \in X^{\otimes m}$ with $\|f\|=1$ such that  
$$
 \|T^{1/2}f \|_2^2 = \|(f|Tf)_A\|  \ge \|T\| -\varepsilon 
$$
Hence we have 
$\|T^{1/2}f \|_2^2 \ge \|E(w)\| - 3 \varepsilon $.
Define $a = S_f^* T S_f = (f|Tf)_A \in A$.  
Then $\|a\| \ge \|E(w)\| -3 \varepsilon  > \varepsilon$. 
By Lemma \ref{lemma:u-epsilon},  there exists
$n \in \mathbb{N}$ and $u \in X^{\otimes n}$ 
sucn that  
$$
 \| u \| _2 \le (\|a\| - \varepsilon)^{-1/2} \qquad \text{and} 
\quad  S_u^* a S_u = I
$$
Then $\|u\| \le (\|E(w)\| -3 \varepsilon)^{-1/2}$. 
The rest of the proof is exactly the same as in \cite{KW} 
Theorem 3.8.  We have 
\[
\| S_f^* PwP S_f -a \|  \le \|S_f\|^2 \|P\|^2 \|w -b\| <\varepsilon 
\]
Therefore  
\[
\|S_u^*S_f^* PwPS_fS_u -I\|  
  < \|u\|^2 \varepsilon 
  \le \varepsilon \|E(w) -3\varepsilon \|^{-1} \le 1.
\]
Hence  
$S_u^*S_f^*PwPS_fS_u$ is invertible. Thus there  exists 
$v \in {\mathcal O}_X$ with $S_u^*S_f^*PwPS_fS_u v =I$. 
Put $z_1=S_u^*S_f^*P$ and $z_2=PS_fS_u v$. Then  
$z_1 w z_2 = I$.  

\end{proof}

\noindent
{\bf Remark}. J. Schweizer \cite{S} showed that ${\mathcal O}_X$ 
is simple  
if a Hilbert bimodule $X$ is 
minimal and non-periodic.  Any 
$X$-invariant ideal $J$ of $A$  corresponds to a closed subset 
$F$ of $K$ with $\sum _i \gamma_i (F) \subset F$. Since such a 
closed set $F$ is $\phi$ or $K$, $X$ is minimal.
Since $A$ is commutative 
and $L(X_A)$ is non-commutative, $X$ is non-periodic.  Thus 
Schweizer's theorem also implies that 
${\mathcal O}_{\gamma}(K)$ is simple. 
Our theorem gives simplicity and pure infiniteness with a 
direct proof.

\begin{prop} Let $(K,d)$ be a compact separable metric space 
and  $\gamma = (\gamma_1,\dots , \gamma_N)$ be a system of proper contractions
on $K$.  Assume that $K$ is self-similar. 
Then the associated $C^*$-algebra 
${\mathcal O}_{\gamma}(K)$ 
is separable and nuclear, and satisfies the Universal  
Coefficient Theorem.  
\end{prop}

\begin{proof}Since ${\mathcal J}_X$ and  ${\mathcal T}_X$ are 
KK-equivalent to abelian $C^*$-algebras $I_X$ and $A$, 
the quotient ${\mathcal O}_X \cong {\mathcal T}_X/{\mathcal J}_X$ 
satisfies the UCT. Also ${\mathcal O}_X$ is shown to be nuclear 
as in an argument of \cite{DS}.   
\end{proof}

\noindent
{\bf Remark}. In the above situation the isomorphisms class of
${\mathcal O}_{\gamma}(K)$
is completely determined by the $K$-theory together with the class of the 
unit by the classification theorem by Kirchberg-Phillips \cite{Ki},
\cite{Ph}.

\section{examples}
We collect  typical  examples from a fractal geomery. 
We also give a general condition that 
the associated $C^*$-algebra 
${\mathcal O}_{\gamma}(K)$ 
is isomorphic to a 
Cuntz algebra ${\mathcal O}_N$.   

We shall calculate the K-groups  by the following six-term exact 
sequence due to Pimsner \cite{Pi}.
$$
\begin{CD}
   K_0(I_X) @>{id - [X]}>> K_0(A) @>i_*>> 
K_0({\mathcal O}_{\gamma}(K)) \\
   @A{\delta _1}AA
    @.
     @VV{\delta _0}V \\
   K_1({\mathcal O}_{\gamma}(K)) @<<i_*< K_1(A) 
@<<{id - [X]}< K_1(I_X)
\end{CD}
$$

\noindent
{\bf Example 4.1.} (Cantor set). 
Let $\Omega = [0,1]$ and $\gamma_1$ and $\gamma_2$ be two 
contractions defined by 
\[
\gamma_1(y) = \frac{1}{3}y  \text{  and  } 
\gamma_2(y) = \frac{1}{3}y + \frac{2}{3}.
\]
Then the self-similar set $K = K(\gamma_1,\gamma_2)$ is the 
Cantor set and the associated $C^*$-algebra 
${\mathcal O}_{(\gamma_1,\gamma_2)}(K)$ is isomorphic to a 
Cuntz algebra ${\mathcal O}_2$.

\medskip

\noindent
{\bf Example 4.2.} (full shift) 
The full $N$-shift space $\{1,2,...,N\}^{\mathbb N}$ is the  space 
of one-sided sequences $x = (x_n)_{n \in {\mathbb N}}$ of symbols
 $\{1,2,...,N\}$. 
Define a system $\sigma = (\sigma_1, \dots , \sigma_N)$ of 
$N$ contractions  on 
$\{1,2,...,N\}^{\mathbb N}$ by 
$$
\sigma_j(x_1,x_2,..., ) = (j,x_1,x_2,..., )
$$
Then each $\sigma_j$ is a proper contraction with the Lipschitz 
constant $Lip(\sigma _j) = \frac{1}{2}$.  
The self-similar set  $K(\sigma_{1},\sigma_{2},..,\sigma_{N}) = 
\{1,2,...,N\}^{\mathbb N}$.  
Then associated $C^*$-algebra 
${\mathcal O}_{\sigma}(K)$ 
is isomorphic to a 
Cuntz algebra ${\mathcal O}_N$  as in \cite{PWY} section 4.

\medskip

\noindent
{\bf Definition.} Recall that a system 
$\gamma = (\gamma_1,\dots , \gamma_N)$ satisfies the 
{\it strong separation condition} in $K$ 
if 
\[
K = \cup _{i=1}^N \gamma (K)   \text{ and } 
\gamma _i (K) \cap \gamma _j(K) = \phi 
\text{ for } i \not= j .
\]  

We say that a system 
$\gamma = (\gamma_1,\dots , \gamma_N)$ satisfies the 
{\it graph separation condition} in $K$ 
if 
\[
K = \cup _{i=1}^N \gamma (K)   \text{ and } 
cograph \ \gamma _i\cap cograph \ \gamma _j = \phi 
\text{ for } i \not= j ,
\] 
where  
$cograph \ \gamma _i := \{(x,y) \in K^2 ; x = \gamma _i(y)\}$. 
It is clear that \\
(strong separation condition) $\Rightarrow$ (graph separation condition) and \\ 
(strong separation condition)$\Rightarrow$ (open set condition), but 
the converses are not true in general. 

If a system 
$\gamma = (\gamma_1,\dots , \gamma_N)$ satisfies the 
strong separation condition in $K$, then 
the map 
$\pi : \{1,2,...,N\}^{\mathbb N} \rightarrow K$ defined by 
$\{\pi (x)\} =  \cap _{m \geq 1} K_{(x_1, \dots ,x_m)}$ 
is a homeomorphism.   
Since 
$\pi \circ \sigma _i = \gamma _i \circ \pi$ for $i = 1, \dots , N$, 
we can identify the system 
$\gamma = (\gamma_1,\dots , \gamma_N)$ with  the system of 
system $\{\sigma_j:j=1,2,..,N \}$ in Example 4.2 (full shift). 
Therefore it is trivial that the  $C^*$-algebra 
${\mathcal O}_{\gamma}(K)$ is  
isomorphic to a 
Cuntz algebra ${\mathcal O}_N$.

\begin{prop} Let $(K,d)$ be a compact metric space 
and  $\gamma = (\gamma_1,\dots , \gamma_N)$ be a system of proper contractions
on $K$.  Assume that $K$ is self-similar. 
If a system 
$\gamma = (\gamma_1,\dots , \gamma_N)$ satisfies the 
graph separation condition,  
then the associated $C^*$-algebra 
${\mathcal O}_{\gamma}(K)$ 
is  isomorphic to a 
Cuntz algebra ${\mathcal O}_N$.
\end{prop}

\begin{proof}
Let $\ _{\beta _i}A$ be a Hilbert bimodule over $A$ 
defined as $\ _{\beta _i}A = A$ as a vector space, 
$a\cdot f \cdot b= \beta _i(a)fb$
for $a,b \in A$ and $f \in \ _{\beta _i}A$ and 
a $A$-valued inner product is the usual one: 
$(f|g)_A = f^*g$
for $f,g \in \ _{\beta _i}A$.  

Let ${\mathcal G}_i := cograph \ \gamma _i$.  
Then $C({\mathcal G}_i)$  a Hilbert bimodule over $A$ by  
\[
(a\cdot f_i \cdot b)(\gamma_i(y),y) 
= a(\gamma_i(y))f(\gamma_i(y),y)b(y)
\]
for $a,b \in A$ and $f_i \in C({\mathcal G}_i)$. 
An $A$-valued 
inner product $(\ |\ )_A$ is defined  by 
\[
(f_i|g_i)_A(y) =  
\overline{f(\gamma _i(y),y)}g(\gamma _i(y),y) 
\]
for $f_i,g_i \in C({\mathcal G}_i)$ and $y \in K$.
It is clear that there exists a $A$-$A$ bimodule isomorphism 
$\psi : \ _{\beta _i}A \rightarrow C({\mathcal G}_i)$
preserving $A$-valued inner product such that 
$\psi(f)(\gamma _i(y),y) = f(y)$ for $f \in \ _{\beta _i}A$
and $y \in K$. 
Since the system 
$\{\gamma_j:j=1,2,..,N \}$ satisfies the 
graph separation condition, we have isomorphisms 
\[
C({\mathcal G}) \cong \oplus _{i=1}^N C({\mathcal G}_i) \cong 
 \oplus _{i=1}^N \ _{\beta _i}A .
\] 
Since each $\gamma _i$ is a proper contraction, 
$C^*$-algebra 
${\mathcal O}_{\gamma}(K)$ 
is  isomorphic to a 
Cuntz algebra ${\mathcal O}_N$ by \cite{PWY} section 4. 
\end{proof}

\noindent
{\bf Example 4.3.} (branches of the inverse of a tent map) 
A tent map $h : [0,1] \rightarrow [0,1]$ is defined by 
\[
h(x) = \begin{cases}
 2x, & \qquad 0 \leq x \leq \frac{1}{2}, \\
-2x + 2, & \qquad \frac{1}{2} \leq x \leq 1.
\end{cases}
\]
Let 
\[
\gamma_1(y) = \frac{1}{2}y  \text{ and } 
\gamma_2(y) = -\frac{1}{2}y +1 .
\]
Then $\gamma_1$ and $\gamma_2$ are branches of $h^{-1}$. 
The self-similar set  $K(\gamma_{1},\gamma_{2}) = [0,1]$
And the associated $C^*$-algebra 
${\mathcal O}_{(\gamma_{1}, \gamma_{2})}(K)$ 
is isomorphic to the $C^*$-algebra ${\mathcal O}_{z^2-2}$ 
associated a polynomial $z^2-2$ \cite{KW}. Since the K-groups 
$(K_0, K_1, [1])$ with the posiiton of the unit $[1]$ in $K_0$ 
are equal, it is also isomorphic to 
a Cuntz algebra ${\mathcal O}_{\infty}$. The 
system $(\gamma_1, \gamma_2)$ satisfies   
open set condition but does not satisfies 
graph separation condition. 

We modify the example a bit.  Let 
\[
\gamma_1'(y) = \frac{1}{2}y  \text{ and } 
\gamma_2'(y) = \frac{1}{2}y +\frac{1}{2} .
\]
Then $\gamma_1'$ and $\gamma_2'$ are not branches of the 
inverse of a certain function, because 
$\gamma_1'(1) = \gamma_2'(0) = \frac{1}{2}$.  
The self-similar set  $K(\gamma_1',\gamma_2') = [0,1]$.
The system $(\gamma_1' \gamma_2')$ satisfies   
graph separation condition but does not satisfy 
strong separation condition.  
And the associated $C^*$-algebra 
${\mathcal O}_{(\gamma_1',\gamma_2')}(K)$ 
is isomorphic to the Cuntz algebra ${\mathcal O}_2$ .

\medskip

\noindent
{\bf Example 4.4.}(Koch curve) Let 
$\omega = \frac{1}{2} + i \frac{\sqrt{3}}{6} \in {\mathbb C}$. 
Consider two contractions 
$\gamma _1, \gamma _2$ on the triangle domain 
$\triangle \subset {\mathbb C}$ with 
vertices $\{ 0, \omega, 1\}$  defined by 
$\gamma _1(z) = \omega \overline{z}$ and 
$\gamma _2(z) = (1-\omega)(\overline{z} - 1) +1$, 
for $z \in {\mathbb C}$.  Then the 
self-similar set $K$ is called the Koch curve.  But 
these two contractions are not inverse branches of a map 
on $K$ 
because $\gamma_1(1) = \gamma_2(0) = \omega$. 
We modify the 
construction of contractions.  Put  
$\tilde{\gamma _1} = \gamma _1$, 
$\tilde{\gamma _2} = \tau \circ \gamma _2$,  
where $\tau$ is to turn over. 
Then $\tilde{\gamma _1}, \tilde{\gamma _2}$ 
are inverse branches of a map $h$ on K. 
$C^*$-algebra ${\mathcal O}_{(\gamma _1,\gamma _2)}(K)$
is isomorphic to the Cuntz algebra 
${\mathcal O}_2$. Moreover $C^*$-algebra 
${\mathcal O}_{(\gamma_w)_{w \in W_n}}(K)$
is isomorphic to the Cuntz algebra 
${\mathcal O}_{2^n}$. Moreover 
$C^*$-algebra 
${\mathcal O}_{(\tilde{\gamma}_w)_{w \in W_n}}(K)$
is isomorphic to  a purely infinite, simple 
$C^*$-algebra ${\mathcal O}_{T_{2^n}}([0,1])$, where 
$T_n$ is the Tchebychev polynomials defined by 
$\cos nz = T_n(\cos z)$, see Example 4.5. in \cite{KW}.
Thus we  have 
$K_0({\mathcal O}_{(\tilde{\gamma}_w)_{w \in W_n}}(K)) = {\mathbb Z}^{2^n-1}$ and 
$K_1({\mathcal O}_{(\tilde{\gamma}_w)_{w \in W_n}}(K))) = 0$. 

\medskip

\noindent
{\bf Example 4.5.}(Sierpinski gasket) 
Recall that the usual Sierpinski 
gasket $K$ is constructed by three contractions 
$\gamma _1, \gamma _2, \gamma _3$ on the regular triangle $T$
in ${\mathbb R}^2$ with three vertices $P = (1/2,\sqrt{3}/2)$, 
$Q = (0,0)$ and $R = (1,0)$ such that 
$\gamma _1(x,y) = (\frac{x}{2} + \frac{1}{4}, \frac{y}{2} + \frac{\sqrt{3}}{4})$, $\gamma _2(x,y) = (\frac{x}{2}, \frac{y}{2})$, 
$\gamma _3(x,y) = (\frac{x}{2} + \frac{1}{2}, \frac{y}{2})$. 
Then the  self-similar set $K$ is called 
a Sierpinski gasket.  But these three contractions are not 
inverse branches of a map, because 
$\gamma _1(Q) = \gamma _2(P)$. 

Ushiki \cite{U} discovered a rational 
function whose Julia set is homeomorphic to the  
Sierpinski gasket.  See also \cite{Kam}.  For example,   
let $R(z) = \frac{z^3-\frac{16}{27}}{z}$. Then the 
Julia set $J_R$ is 
homeomorphic to the  Sierpinski gasket K and $J_R$ contains 
three critical points. 
Therefore we need to modify the 
construction of contractions.  Put  
$\tilde{\gamma _1} = \gamma _1$, 
$\tilde{\gamma _2} = \alpha _{-\frac{2\pi}{3}} \circ \gamma _2$,  
and $\tilde{\gamma _3} = \alpha _{\frac{2\pi}{3}} \circ \gamma _3$, 
where $\alpha _{\theta}$ is a rotation by the angle $\theta$. 
Then $\tilde{\gamma _1}, \tilde{\gamma _2},  \tilde{\gamma _3}$ 
are inverse branches of a map $h: K \rightarrow K$, which is 
conjugate to $R : J_R \rightarrow J_R$.  
Then $C^*$-algebra ${\mathcal O}_R \cong 
{\mathcal O}_{(\tilde{\gamma _1}, \tilde{\gamma _2},  \tilde{\gamma _3})}(K)$ 
is  a purely infinite, simple $C^*$-algebra, 
and $K_0({\mathcal O}_R)$ contains a torsion free element.  But 
$C^*$-algebra ${\mathcal O}_{(\gamma _1,\gamma _2, \gamma _3)}(K)$.
is isomorphic to the Cuntz algebra 
${\mathcal O}_3$, because 
the system $(\gamma _1,\gamma _2, \gamma _3)$  satisfies   
graph separation condition. 
 Therefore  $C^*$-algebra
 ${\mathcal O}_{(\gamma _1,\gamma _2, \gamma _3)}(K)$ and 
${\mathcal O}_{(\tilde{\gamma _1}, \tilde{\gamma _2},  \tilde{\gamma _3})}(K)$ are not isomorphic.  See \cite{KW}. 

\medskip

\noindent
{\bf Example 4.6.}(Sierpinski carpet) 
Recall that the usual Sierpinski 
carpet $K$ is constructed by eight contractions 
$\gamma _1,\dots , \gamma _8$ on the regular square 
$S = [0,1]\times [0,1]$
in ${\mathbb R}^2$ with four vertices $P_1 = (0,1)$, 
$P_2 = (0,0)$, $P_3 = (1,0)$ and $P_4 = (1,1)$ such that  
$\gamma _1(x,y) = (\frac{x}{3}, \frac{y}{3})$,
$\gamma _2(x,y) = (\frac{x}{3} + \frac{1}{3}, \frac{y}{3})$, 
$\gamma _3(x,y) = (\frac{x}{3} + \frac{2}{3}, \frac{y}{3})$, 
$\gamma _4(x,y) = (\frac{x}{3}, \frac{y}{3} + \frac{1}{3})$,
$\gamma _5(x,y) = (\frac{x}{3} + \frac{2}{3}, \frac{y}{3} + \frac{1}{3})$,
$\gamma _6(x,y) = (\frac{x}{3}, \frac{y}{3} + \frac{2}{3})$,
$\gamma _7(x,y) = (\frac{x}{3} + \frac{1}{3}, \frac{y}{3} + \frac{2}{3})$,
$\gamma _8(x,y) = (\frac{x}{3} + \frac{2}{3}, \frac{y}{3} + \frac{2}{3})$.
Then the  self-similar set $K$ is called 
a Sierpinski carpet.  But these eight contractions are not 
continuous branches of the inverse of any map $h : K \rightarrow K$, 
because $\gamma _1(P_1) = \gamma _4(P_2)$. 
We shall  modify the 
construction of contractions as follows:
$\gamma'_1(x,y) = \gamma_1(x,y)$,
$\gamma'_2(x,y) = (-\frac{x}{3} + \frac{2}{3}, \frac{y}{3})$, 
$\gamma'_3(x,y) = \gamma_3(x,y)$, 
$\gamma'_4(x,y) = (\frac{x}{3}, -\frac{y}{3} + \frac{2}{3})$,
$\gamma'_5(x,y) = (\frac{x}{3} + \frac{2}{3}, -\frac{y}{3} + \frac{2}{3})$,
$\gamma'_6(x,y) = \gamma_6(x,y)$,
$\gamma'_7(x,y) = (-\frac{x}{3} + \frac{2}{3}, \frac{y}{3} + \frac{2}{3})$,
$\gamma'_8(x,y) = \gamma_8(x,y)$.
Then their  self-similar set is the same 
Sierpinski carpet $K$ as above. 
And $\gamma'_1, \dots, \gamma'_8$ are  
continuous branches of the inverse of a map $h : K \rightarrow K$. 
Since 
\begin{align*}
B &= B(\gamma'_1, \dots, \gamma'_8) \\
&=  (([0,\frac{1}{3}] \cup [\frac{2}{3},1])\times 
\{\frac{1}{3},\frac{2}{3}\})
\cup (\{\frac{1}{3},\frac{2}{3}\}\times 
([0,\frac{1}{3}] \cup [\frac{2}{3},1])), 
\end{align*}
Hence $K_0(C(B)) \cong {\mathbb Z}^4$ and  $K_1(C(B)) \cong 0$.
Since we have $K_0(C(K)) \cong {\mathbb Z}$ and 
$K_1(C(K)) \cong {\mathbb Z}^{\infty}$, 
$K_0({\mathcal O}_{(\gamma'_1, \dots, \gamma'_8)}(K))$ 
contains a torsion free element.
But $C^*$-algebra ${\mathcal O}_{(\gamma _1,\dots, \gamma _8)}(K)$.
is isomorphic to the Cuntz algebra ${\mathcal O}_8$, 
because the system $(\gamma _1,\dots, \gamma _8)$  satisfies   
graph separation condition. 
 Therefore  purely infinite, simple $C^*$-algebras 
${\mathcal O}_{(\gamma _1,\dots, \gamma _8)}(K)$ and 
${\mathcal O}_{(\gamma'_1, \dots, \gamma'_8)}(K)$
are not isomorphic.


%

\end{document}